\documentclass[12pt,a4paper]{amsart}

\usepackage{tikz}
\tikzset{
op/.style = {shape=rectangle, rounded corners, draw=black, align=center,
fill=blue!30},base align tikzpicture/.style={execute at end picture={
\path (current bounding box.north) -- (current bounding box.south)
node[midway](Xphantom){\phantom{X}};},baseline={(Xphantom.base)}}
        }
\usepackage{tikz-cd}

\usepackage{graphicx} 
\usepackage{hyperref}
\usepackage[notref,notcite]{}
\hypersetup{
    colorlinks=true,
    linkcolor=blue,
    filecolor=magenta,     
    urlcolor=cyan,
    pdftitle={Codimbraid},
    }
\usepackage{amssymb,amsmath,amsthm,amsfonts,mathrsfs,bbm}
\usepackage{array,amssymb,amsxtra} 
\usepackage{geometry}
\usepackage[shortlabels]{enumitem} 
\usepackage[normalem,normalbf]{ulem}
\usepackage[all]{xy} 
\usepackage[american]{babel}
\usepackage{fancybox}
\usepackage{diagbox}
\usepackage{thmtools}
\usepackage{cite}
\usepackage{yhmath}
\newcommand{\N}{\mathbb{N}}

\newcommand{\OO}{\mathscr{O}_{p, q}(\normalfont\text{SL}_{2}(k))}

\newcommand{\Naksp}{\mathrm{NSp}}
\usepackage{thmtools}

\DeclareMathOperator{\End}{End}

\newtheorem{theorem}{Theorem }[section]
\newcommand{\bt}{\begin{theorem}}
\newcommand{\et}{\end{theorem}}

\definecolor{navajowhite}         {rgb}{1.00,0.87,0.68}
\newtheorem{prop}[theorem]{Proposition}
\newcommand{\bp}{\begin{prop}}
\newcommand{\ep}{\end{prop}}

\newtheorem{dep}[theorem]{Proposition-Definition}
\newcommand{\bdep}{\begin{dep}}
\newcommand{\edep}{\end{dep}}

\newtheorem{cor}[theorem]{Corollary}
\newcommand{\bc}{\begin{cor}}
\newcommand{\ec}{\end{cor}}
\definecolor{babyblueeyes}{rgb}    {0.63, 0.79, 0.95}
\newtheorem{lemma}[theorem]{Lemma}
\newcommand{\bl}{\begin{lemma}}
\newcommand{\el}{\end{lemma}}

\theoremstyle{definition}

\newtheorem{dfn}[theorem]{Definition}
\newcommand{\bd}{\begin{dfn}}
\newcommand{\ed}{\end{dfn}}

\theoremstyle{remark}

\newtheorem{rmq}[theorem]{Remark}
\newcommand{\brq}{\begin{rmq}}
\newcommand{\erq}{\end{rmq}}

\theoremstyle{definition}

\newtheorem{ex}[theorem]{Example}
\newcommand{\bex}{\begin{ex}}
\newcommand{\eex}{\end{ex}}

\newcommand{\bpf}{\begin{proof}}
\newcommand{\epf}{\end{proof}}

\title[Coordinate algebras of the braided quantum groups SL$_2$]{Isomorphic classification of the coordinate algebras of the two-parameter braided quantum groups SL$_2$}
\author{Julien Bichon}
\address{Universit\'e Clermont Auvergne, CNRS, LMBP, F-63000 CLERMONT-FERRAND, FRANCE}
\email{julien.bichon@uca.fr}
\author{Thi Hoa Emilie Nguyen}
\email{thi\_hoa.nguyen@uca.fr}

\thanks{Work funded by ANR-25-CE40-5010}

\begin{document}

\begin{abstract}
We classify, up to isomorphism, the coordinate algebras of the two-parameter braided quantum groups $\textrm{SL}_2$.
\end{abstract}

\maketitle

\section{Introduction}
We work over a fixed base field $k$, and for $p,q\in k^*$, we consider the algebra
$\OO$  presented by generators $a, b, c, d$ subject to the relations
\[ba = qab, \ ca = pac, \ db = qbd, \ dc = pcd, \ bc=cb, \ ad - p^{-1}bc = 1 = da - qbc.\]

This algebra was considered in our previous paper \cite{BN26}. 
For $p=q$,  $\mathscr{O}_{q, q}(\normalfont\text{SL}_{2}(k))=\mathscr{O}_{q}(\normalfont\text{SL}_{2}(k))$ has a natural Hopf algebra structure, and is the familiar coordinate algebra of the quatum group $\text{SL}_2$. For $p\not=q$, it was shown in\cite{BN26}  that $\OO$ carries  a Hopf algebra structure in the braided category of $\mathbb Z$-graded vector spaces (with  braiding associated with the scalar $\xi=p^{-1}q$), and hence this braided Hopf algebra was called there
the coordinate algebra of the two-parameter braided quantum group $\text{SL}_2$. 

The algebras $\OO$ were studied in \cite{BN26}, with this study culminating in showing that these are twisted Calabi-Yau algebras of dimension $3$. However, while we were aware at the time of writing \cite{BN26} that there are some $\OO$ that cannot be isomorphic to a single parameter version,   the  isomorphic classification was not discussed at all there. The aim of the present paper is to fill this gap, as follows.

\begin{theorem}\label{thm:iso}
 For $p,q, \alpha, \beta \in k^*$, the algebras $\OO$ and $\mathscr{O}_{\alpha, \beta}(\normalfont\text{SL}_{2}(k))$ are isomorphic if and only if $(\alpha, \beta)\in \{(p,q), (q,p), (p^{-1},q^{-1}), (q^{-1},p^{-1})\}$.  
\end{theorem}


We also take the opportunity of this paper to introduce what seems to be a  useful invariant for twisted Calabi-Yau algebras having a central group of invertible elements: the Nakayama spectrum, that is the set of eigenvalues of the Nakayama automorphism. In the case of $\OO$, the Nakayama spectrum is the subgroup $\langle pq\rangle\subset k^*$. This complements well the Alev-Dumas invariant \cite{AD94}, which, in the $\OO$ case,  is the subgroup $\langle p, q \rangle \subset k^*$.

This paper is organized as follows. In Section \ref{sec:basic} we recall and complement some basic results about $\OO$. In Section 
\ref{sec:invariants} we discuss two general invariants for algebras: we recall the Alev-Dumas invariant,  and we introduce the Nakayama spectrum. Both these invariants are computed for $\OO$. In the final Section \ref{sec:isoclass}, we discuss the conormal bimodule of an ideal, explain how it can be used to produce an invariant via the spectrum of certain automorphisms, and complete the proof of Theorem \ref{thm:iso}.


\subsection*{Acknowledgements.} In this paper we collaborated with M365 Copilot. The collaboration started by asking Copilot to search for  a proof of the $p=q$ case of Theorem \ref{thm:iso} in the literature, a folk-known result for which, strangely, we still do not know a printed reference. Copilot first furnished proofs in the non root of unity case (one was very computational, the second one essentially amounts to what we present in Subsection \ref{subsec:AD}, and that we were already aware of), and after a few other requests, furnished a proof working as well at roots of unity. We checked the proof, which appeared to be correct although slightly clumsy. We then polished and reworked  these arguments in order that they could be adapted to  the $p\not=q$ case, yielding Lemma \ref{lem:conormal} as a result, that we then used in a crucial way to prove finally Theorem \ref{thm:iso} (the $p\not= q$ case required several subcases with additional complications). 

On the Human side, it is a pleasure to thank François Dumas for his interest in this project and for  discussions on these questions.

\section{Basic facts about $\OO$}\label{sec:basic}

The aim of  this section is to provide some useful facts about $\OO$. The first one is an explicit vector space basis of $\OO$, obtained, as usual, using the diamond lemma.

\bp\label{prop:diamond2}
The set $\{a^{i}b^{j}c^{k} \mid i, j, k \in \N\} \cup \{b^{j}c^{k}d^{l} \mid j, k \in \N, l \in \N^{*}\}$ is a vector space basis of $\OO$.
\ep


\bp \label{domain}
The  algebra $\OO$  is a Noetherian domain.
\ep

\bpf
The result is given in \cite{BN26}. We provide, for the sake of completeness, the proof, which follows classical arguments as in \cite[I.1]{BG02}. We consider first the algebra $A_{p,q}$ presented by generators 
 $a, b, c, d$ and relations
\[ba = qab, \,ca = pac, \,bc = cb, \,
db = qbd, \,dc = pcd, \,ad - da = (p^{-1} - q)bc\]
and remark that $A_{p,q}$ is an iterated Ore extension, hence is a Noetherian domain.
Put  $D_{p, q} = ad - p^{-1}bc = da - qbc$. We have $D_{p, q} \in Z(A_{p,q})$ and  
$\OO = A_{p,q} /(D_{p, q} - 1)$, hence $\OO$ is Noetherian. Consider then the localization $A_{p,q}[D^{-1}_{p,q}]$ with respect to the central regular element $D_{p, q}$. We then have an algebra isomorphism
\begin{align*}
f: A_{p,q}[D_{p, q}^{-1  }] &\longrightarrow \OO \otimes k[z, z^{-1}]\\
a,b,c,d &\longmapsto a\otimes z, b\otimes z, c\otimes 1, d\otimes 1
\end{align*}
and since $A_{p,q}[D^{-1}_{p,q}]$ is a domain, so is $\OO$.
\epf

The easy part of Theorem \ref{thm:iso} is dealt with in the following result.

\bp \label{prop:isobasic}  For any $p,q\in k^*$,
we have algebra isomorphisms \[\OO \simeq \mathscr{O}_{q, p}(\normalfont\text{SL}_{2}(k)) \simeq \mathscr{O}_{p^{-1}, q^{-1}}(\normalfont\text{SL}_{2}(k)) \simeq \mathscr{O}_{q^{-1}, p^{-1}}(\normalfont\text{SL}_{2}(k))\]
\ep

\begin{proof}
 Write $\mathscr{O}_{p,q} =\OO$. The announced isomorphisms are given as follows: 
 \begin{align*}
\mathscr{O}_{p,q} &\longrightarrow \mathscr{O}_{q,p}  \quad & \mathscr{O}_{p,q} &\longrightarrow \mathscr{O}_{q^{-1},p^{-1}} \quad & \mathscr{O}_{p,q} &\longrightarrow \mathscr{O}_{p^{-1},q^{-1}}\\
\begin{pmatrix}
a &b\\
c &d
\end{pmatrix}
&\longmapsto
\begin{pmatrix}
a & pq^{-1}c \\
b & d
\end{pmatrix}
&\begin{pmatrix}
a &b\\
c &d
\end{pmatrix}
&\longmapsto \begin{pmatrix}
d & c \\
b & a
\end{pmatrix} 
& \begin{pmatrix}
a & b \\
c & d
\end{pmatrix} 
& \longmapsto
\begin{pmatrix}
d & pq^{-1}b \\
c & a
\end{pmatrix}
\end{align*}
The details are left to the reader.
\end{proof}

Our next result describes the abelianization of $\OO$, and provides
  a first very rough classification, to be used in Section \ref{sec:isoclass}. The details of the proof are left to the reader.

\bp \label{prop:abelianization}
For $p,q \in k^*$, put $A=\OO$. Let $I$ denote the commutator ideal of $A$ and  $A_{\rm ab}$  the abelianization of $A$.
\begin{enumerate}[(a)]
 \item For $p\not=1$, $q\not=1$ and $pq\not=1$, we have $I=Ab+Ac$ and $A_{\rm ab}\simeq k[t^{\pm 1}]$ with the above isomorphism induced by the algebra map 
 \begin{align*}
  A &\longrightarrow k[t^{\pm 1}] \\ 
  \begin{pmatrix}
a &b\\
c &d
\end{pmatrix} & \longmapsto
\begin{pmatrix}
t &0\\
0 &t^{-1}
\end{pmatrix}
 \end{align*}
\item For $p=1$ and $q\not=1$, we have $I=Ab$ and $A_{\rm ab}\simeq k[t^{\pm 1}, y]$ with this isomorphism induced by the algebra map 
 \begin{align*}
  A &\longrightarrow k[t^{\pm 1}, y] \\ 
  \begin{pmatrix}
a &b\\
c &d
\end{pmatrix} & \longmapsto
\begin{pmatrix}
t &0\\
y &t^{-1}
\end{pmatrix}
 \end{align*}
 \item  For $q=p^{-1}\not=1$, we have $I=Aab + Aac+ Abd+ Acd$ and $A_{\rm ab}\simeq k[t^{\pm 1}]\times k[t^{\pm 1}]$ with this isomorphism induced by the algebra map 
 \begin{align*}
  A &\longrightarrow k[t^{\pm 1}]\times k[t^{\pm 1}] \\ 
  \begin{pmatrix}
a &b\\
c &d
\end{pmatrix} & \longmapsto
\begin{pmatrix}
(t,0) & (0,t)\\
(0,-pt^{-1})  & (t^{-1},0)
\end{pmatrix}
\end{align*}
\end{enumerate}
\ep

The group of invertible elements of an algebra $A$ is denoted $A^\times$. The rest of this section is devoted to compute $\OO^\times$, which will be needed in Section \ref{sec:invariants}. 
For this we first consider the two-parameter quantum matrix algebra $\mathscr{O}_{q , p}(M_2(k))$, see \cite{Tak90}, presented by generators $a,b,c,d$ and relations
\[ba = qab, \ ca = pac, \ db = pbd, \ dc = qcd, \ pbc=qcb, \
da-ad = (p-q^{-1})bc.\]
The algebra $\mathscr{O}_{q , p}(M_2(k))$ is a domain since it is an iterated Ore extension of domains, and since it is graded connected as well, we have $\mathscr{O}_{q , p}(M_2(k))^\times =k^*$.  
The quantum determinant is defined by 
		\[\det = ad - q^{-1}bc = da - pbc\]
We have
		\[a\det = \det a, \ b\det = p^{-1}q\det b, \
		c\det = pq^{-1}\det c, \ d\det = \det d\]
		so that $\det$ is a normal element in $\mathscr{O}_{q , p}({M}_2(k))$.
We consider the algebra $\mathscr{O}_{q , p}({M}_2(k)) /(\det)$

\bl
 The set $\{a^{i}b^{j}c^{k} \mid i, j, k \in \N\} \cup \{b^{j}c^{k}d^{l} \mid j, k \in \N, l \in \N^{*}\}$
  is a vector space basis of $\mathscr{O}_{q , p}({M}_2(k)) /(\det)$. Moreover  $\mathscr{O}_{q , p}({M}_2(k)) /(\det)$ is a domain.		
\el

\bpf
The announced basis is obtained by using the diamond lemma. To prove that $\mathscr{O}_{q , p}({M}_2(k)) /(\det)$ is a domain, we follow closely  the argument of Goodearl and Lenagan \cite[Theorem 1.1]{GL00} in the case $p=q$.
It is straightforward to check that there exists a unique morphism of algebras
\begin{align*}
			\theta : \mathscr{O}_{q , p}({M}_2(k)) &\longrightarrow k_{p}[x_{1}, x_{2}] \otimes k_{q}[y_{1}, y_{2}]\\
			\begin{pmatrix}
				a & b\\
				c & d
			\end{pmatrix} &\longmapsto \begin{pmatrix}
				x_{1} \otimes y_{1}& x_{1} \otimes y_{2}\\
				x_{2} \otimes y_{1} & x_{2} \otimes y_{2}
			\end{pmatrix}
		\end{align*}
		We have  $\theta(ad - q^{-1}bc) = (x_{1} \otimes y_{1})(x_{2} \otimes y_{2}) - q^{-1}(x_{1} \otimes y_{2})(x_{2} \otimes y_{1}) = 0$. Thus $\det \in \ker(\theta)$ and there exists a unique morphism 
		\begin{align*}
			\bar{\theta} : \mathscr{O}_{q , p}({M}_2(k)) / (\det) &\longrightarrow k_{p}[x_{1}, x_{2}] \otimes k_{q}[y_{1}, y_{2}]\\
			\begin{pmatrix}
				\bar{a} & \bar{b}\\
				\bar{c} & \bar{d}
			\end{pmatrix} &\longmapsto \begin{pmatrix}
				x_{1} \otimes y_{1}& x_{1} \otimes y_{2}\\
				x_{2} \otimes y_{1} & x_{2} \otimes y_{2}
			\end{pmatrix}
		\end{align*}
		Using the respective bases of the algebras, we see that this morphism is injective. Furthermore, $k_{p}[x_{1}, x_{2}] \otimes k_{q}[y_{1}, y_{2}]$ is a domain, since it is isomorphic to a coordinate algebra of a multiparametric quantum affine space, and hence we conclude that $\mathscr{O}_{q , p}({M}_2(k)) / (\det)$ is also a domain. 
\epf
		
		We have seen that $\det$ is a normal element in the Noetherian domain $\mathscr{O}_{q , p}({M}_2(k))$.
		Hence we form the localization $\mathscr{O}_{q , p}({\rm GL}_2(k)) = \mathscr{O}_{q , p}({M}_2(k))[\det^{-1}]$ \cite{Tak90}, which can be presented in terms of generators and relations in the usual manner.
		
\begin{prop}\label{prop:invOOG}
 We have  $\mathscr{O}_{q , p}({\rm GL}_2(k))^\times = \{\lambda \det^n, \ \lambda \in k^*, n \in \mathbb Z\}$
\end{prop}		

\bpf
Let $x \in \mathscr{O}_{q , p}({\rm GL}_2(k))^\times $. We write $x = u\det^m$, for $m \in \mathbb{Z}$ and $u\in \mathscr{O}_{q , p}({M}_2(k))$ such that $u \not \in (\det)$.  Similarly we write $y=x^{-1} = \det^n v $   for $n \in \mathbb{Z}$ and $v\in \mathscr{O}_{q , p}({M}_2(k))$ such that $v \not \in (\det)$.  We then have $yx=1 \Rightarrow vu =\det^{-m-n}$. If $m+n >0$ this contradicts that $u,v \in  \mathscr{O}_{q , p}({M}_2(k))$. If $m+n<0$ we get that $uv \in (\det)$ while $u \not  \in (\det)$ and $v\not \in (\det)$,  which contradicts the fact that  $\mathscr{O}_{q , p}({M}_2(k)) /(\det)$ is a domain. Hence $m+n = 0$ and $uv=1$,  so that  $u$ is invertible in $\mathscr{O}_{q , p}({M}_2(k)) $, and $u=\lambda 1$ for $\lambda \in k^*$.
\epf

\begin{prop}\label{prop:invOO}
 We have $\OO^\times = k^*$.
\end{prop}

\begin{proof}
We have seen in \cite[Proposition 6.6]{BN26} that there is an algebra isomorphism 
 \begin{align*}
 f : \mathscr{O}_{p , q}({\rm GL}_2(k)) & \longrightarrow k\mathbb Z \# \OO \\
\begin{pmatrix}
 a & b \\ c & d
\end{pmatrix}, \ \textrm{det}^{-1}
& \longmapsto 
\begin{pmatrix}
 z\# a & z \# b \\ 1 \# c & 1\# d
\end{pmatrix}, \ z^{-1}\#1
\end{align*}
where the smash product on the right is associated with an appropriate $k\mathbb Z$-module algebra structure on $\OO$ ($z$ is the multplicative generator of $\mathbb Z$). The previous computation of $\mathscr{O}_{q , p}({\rm GL}_2(k))^\times$ thus gives the announced result.
\end{proof}


\section{Two general invariants for algebras}\label{sec:invariants}

In this section we discuss two invariants of general interest for algebras having particular properties. The first one, the Alev-Dumas invariant, is already a well-known tool to prove that some algebras are not isomorphic, and we compute it for $\OO$.  The second one, the Nakayama sprectrum, seems not to have been formally defined  yet, but of course might already be known to some experts.

\subsection{The Alev-Dumas invariant}\label{subsec:AD} The invariant, introduced in \cite{AD94}, is defined as follows.

\bd
Let $A$ be a Noetherian domain, with division ring of fractions  $D={\rm Fract}(A)$. The \textsl{Alev-Dumas invariant} of $A$ is the subgroup of $k^*$ defined by
\[G_{\rm AD}(A) = [D^*,D^*] \cap k^*\]
\ed

It is clear that $G_{\rm AD}(A)$ is an invariant of the Noetherian domain $A$. 

\bex
Let $\mathbf{p} \in M_n(k)$ be a multiplicatively antisymmetric matrix, i.e. 
$p_{ii} = 1=p_{ij}p_{ji}$ for any $i,j$. Consider the  coordinate algebra of the multiparametric quantum affine space
\[\mathscr{O}_{\mathbf{p}}(k^n)= \langle x_1, \ldots , x_n \ | \ 
x_i x_j = p_{ij}x_jx_i,  1\leq i,j \leq n \rangle \]
This is a Notherian domain, and by \cite[Proposition 3.9]{AD94} we have $G_{\rm AD}(\mathscr{O}_{\mathbf{p}}(k^n))= \langle p_{ij}, \ 1 \leq i<j \leq n\rangle$. In particular for the quantum plane algebra $k_q[x,y] = k\langle x,y \ | \ yx = qxy\rangle$ we have $G_{\rm AD}(k_q[x,y])= \langle q \rangle$, so that we can deduce that for $p,q$ non roots of unity, we have $k_q[x,y]\simeq k_p[x,y] \Rightarrow p=q^{\pm 1}$. This is generalized to include the  roots of unity case in \cite{gad14} or in \cite{SAQ15}. As pointed in \cite{gad17}, this becomes quite easy  to prove using \cite[Theorem 1]{BZ17}.
\eex

\bp\label{prop:ADOO}
The Alev-Dumas invariant of $\OO$ is $G_{\mathrm{AD} }(\OO)=\langle p,q \rangle$. Hence
for $p,q, \alpha, \beta \in k^*$, $\OO \simeq \mathscr{O}_{\alpha, \beta}(\normalfont\text{SL}_{2}(k)) \Rightarrow \langle p, q\rangle = \langle \alpha, \beta \rangle$. 
\ep

\bpf Let $D_{p,q} = \operatorname{Fract}(\OO)$. We follow the argument in \cite[Example 1.2.12]{BG02} to identify $D_{p,q}$ with a familiar division ring. Let \[T_{p,q}= \langle x, y, z \mid yx = qxy, zx = pxz, yz = zy \rangle\]
We have $T_{p,q} \simeq \mathscr{O}_\mathbf{p}(k^3)$, for 
 \[\mathbf{p} = \begin{pmatrix}
			1 &q^{-1} &p^{-1}\\
			q &1 &1\\
			p &1 &1
		\end{pmatrix}\]
		In $\OO$, we have
		$
		ad-p^{-1}bc=1,
		$
		and hence
		$
		d=a^{-1}(1+p^{-1}bc)
		$
		in  $D_{p,q}$. Therefore, $D_{p,q}$ is generated, as a division algebra, by $a$, $b$, and $c$. Define an algebra morphism
		\begin{align*}
			f : T_{p, q} &\longrightarrow \OO\\
			x, y, z &\longmapsto a, b, c
		\end{align*}
		Since $T_{p,q}$ is a Noetherian domain and $f$ is injective, $f$ extends uniquely to a homomorphism of division rings
		\[
		\overline{f}:\operatorname{Fract}(T_{p,q}(k))
		\longrightarrow
		\operatorname{Fract}(\OO)=D_{p,q}.
		\]
		This homomorphism is surjective, since $D_{p,q}$ is generated as a division algebra by $a$, $b$, and $c$, and is injective  since any homomorphism between division rings is injective. Hence $\overline{f}$ is an isomorphism.
		We thus have $G_{\textrm{AD}}(T_{p,q}) = G_{\textrm{AD}}(\OO)$, and we conclude by the previous example.
\epf

\brq
It follows from Proposition \ref{prop:ADOO} that, for $p,q\in k^*$ not roots of unity, we have $\mathscr{O}_{p}(\normalfont\text{SL}_{2}(k))\simeq \mathscr{O}_{q}(\normalfont\text{SL}_{2}(k)) \iff p=q^{\pm 1}$. This can also be seen more directly by following the reasoning in \cite{Dan24}. However none of these arguments is sufficient in the root of unity case.
\erq

\brq
It follows from Proposition \ref{prop:ADOO} that if $\langle p,q\rangle$ is a free abelian group of rank two, then $\OO$ is not isomorphic to any single parameter algebra $\mathscr{O}_{\alpha}(\normalfont\text{SL}_{2}(k))$.
\erq

\subsection{The Nakayama spectrum} 
Let $A$ be an algebra. Recall that $A$ is said to be smooth if $A$  admits, as an $A$-bimodule, a finite resolution by finitely generated $A$-bimodules, and that $A$ is said to be twisted Calabi-Yau of dimension $n\geq 0$ if it is smooth and 
 \[H^{i}(A, \,_{A}A \otimes A_{A}) \simeq 
\begin{cases}
	\{0\} \quad \text{if} \quad i \neq n\\
	A_{\mu} \quad \text{if} \quad i = n
\end{cases}
\]
as $A$-bimodules  for an algebra automorphism $\mu \in {\rm Aut}(A)$, called the Nakayama automorphism of $A$. Here $H^*$ denotes Hochschild cohomology, the spaces $H^{i}(A, \,_{A}A \otimes A_{A})$ are naturally $A$-bimodules because $_{A}A \otimes A_{A}$ is an $A^e$-bimodule, and $A_\mu$ has the bimodule structure defined by $a'.a.a'' = a'a\mu(a'')$.

A pleasant feature of twisted Calabi-Yau algebras is that they have a nice Poincaré duality between their homologies and cohomologies \cite{VDB}. We remark below that under an additional assumption, they come equiped with an interesting numerical invariant.

The Nakayama automorphism of a twisted Calabi-Yau algebra is unique up to an inner automorphism, hence is unique if the inner automorphisms all are trivial, i.e.  if the group $A^\times$ of invertible elements of $A$ is contained in the center of $A$. This motivates the following definition.

\bd
Let $A$ be a twisted Calabi-Yau algebra with $A^{\times} \subset Z(A)$.
The \textsl{Nakayama spectrum} of $A$, denoted $\Naksp(A)$, is defined to be the spectrum of the Nakayama automorphism of $A$.
\ed

\bp 
If $A$, $B$ are isomorphic twisted Calabi-Yau algebras with central groups of invertible elements, then $\Naksp(A)= \Naksp(B)$.
\ep

\bpf
The point is to show that
for an algebra isomorphism $\phi : A \to B$, we have $\mu_A= \phi^{-1}\circ \mu_B \circ \phi$, where $\mu_A$, $\mu_B$ are Nakayama automorphisms of $A$, $B$ respectively. The isomorphim $\phi$
induces an isomorphism $A^e \simeq B^e$ between their envelopping algebras, hence an equivalence of categories $\phi^*$ between the corresponding categories of bimodules, and  isomorphisms of $A$-bimodules
\[A_{\phi^{-1}\mu_B \phi} \simeq \phi^*(B_{\mu_B}) \simeq \phi^*({\rm Ext}^n_{B^e}(B,B^e))\simeq {\rm Ext}^n_{A^e}(\phi^*(A), \phi^*(A^e))  \simeq {\rm Ext}^n_{A^e}(A, A^e)\simeq A_{\mu_A}\]
Therefore $\mu_A$ and $\phi^{-1} \circ \mu_B \circ \phi$ only differ by an inner automorphism and hence coincide by the assumption that invertible elements in $A$ are central.
\epf

We know from \cite[Theorem 6.13]{BN26} that $\OO$ is twisted Calabi-Yau of dimension $3$, and from Proposition \ref{prop:invOO} that $\OO^\times = k^*$. Hence we can consider the Nakayama spectrum of $\OO$.

\bp \label{prop:nakspOO}
The Nakayama spectrum of $\OO$ is $\Naksp(\OO)= \langle pq \rangle$. Hence
for $p,q, \alpha, \beta \in k^*$, $\OO \simeq \mathscr{O}_{\alpha, \beta}(\normalfont\text{SL}_{2}(k)) \Rightarrow \langle pq\rangle = \langle \alpha\beta \rangle$. 
\ep

\bpf By  \cite[Theorem 6.13]{BN26}  the Nakayama automorphism of $\OO$ is given by
\begin{align*}
\mu : \quad \OO \quad &\longrightarrow \OO \\
\begin{pmatrix}
a & b\\
c & d
\end{pmatrix} &\longmapsto
\begin{pmatrix}
(pq)^{-1}a &  b\\
c & (pq) d
\end{pmatrix}.
\end{align*}
It thus follows from Proposition \ref{prop:diamond2} that $\Naksp(\OO)= \textrm{Sp}(\mu) = \langle pq \rangle$
\epf

\brq
None of the Alev-Dumas invariant or of the Nakayama spectrum is stronger than the other. Indeed, let $p \in k^*$, $p^2 \not = 1$. The algebras $ \mathscr{O}_{p, p^{-1}}(\normalfont\text{SL}_{2}(k))$ and $ \mathscr{O}_{p}(\normalfont\text{SL}_{2}(k))$ have different Nakayama spectra, while they have the same Alev-Dumas invariant.  On the other hand the algebras $\mathscr{O}_{p^2, 1}(\normalfont\text{SL}_{2}(k))$ and $\mathscr{O}_{p}(\normalfont\text{SL}_{2}(k))$ have the same Nakayama spectrum, but different Alev-Dumas invariants if $p \not \in \langle p^2 \rangle$, i.e. for example if $p$ is not a root of unity.
\erq

\section{Isomorphic classification}\label{sec:isoclass}

\subsection{The conormal bimodule of an ideal}
In this subsection, we explain how, if $A$ is an algebra and $I\subset A$ is a two-sided ideal, the quotient $I/I^2$, the familiar conormal (bi)module from algebraic geometry (see e.g. \cite[Chapter 16]{eis}), can be used to obtain some kind of invariant of $A$, via the spectrum of certain automorphisms. The details are given in the following lemma.

\begin{lemma}\label{lem:conormal}
	Let $A$, $B$ be algebras given with an algebra isomorphism 
	$\phi : A\to B$. Let $I\subset A$ be a (two-sided) ideal and put $J=\phi(I)$. 
	 Consider the $A/I$-bimodule structure on $I/I^2$ (resp. the $B/J$-bimodule structure on $J/J^2$) given, for $a_1,a_2\in A$ and  $x \in I_A$ (resp. $b_1,b_2 \in B$ and $y\in I_B$)  by 
	\[\overline{a_1}\cdot \tilde{x}\cdot \overline{a_2} = \widetilde{a_1xa_2},\quad (\textrm{resp. } \overline{b_1}\cdot \tilde{y}\cdot \overline{b_2} = \widetilde{b_1yb_2})\]
	where $\overline{\cdot}$ and $\widetilde{\cdot}$ stand respectively for the classes of elements of $A$, $B$ in $A/I$, $B/J$, and for the classes of elements of $I$, $J$ in $I/I^2$, $J/J^2$. 
If $t =\overline{a}\in A/I$ is invertible, so is $s=\overline{\phi(a)}\in B/J$, and the linear maps ${\rm ad}(t) \in \End_k(I/I^2)$   and 	${\rm ad}(s) \in \End_k(J/J^2)$, respectively defined by
\[ {\rm ad}(t)(\widetilde{x}) = t\cdot \widetilde{x}\cdot  t^{-1}, \ {\rm ad}(s)(\widetilde{y}) = s\cdot \widetilde{y}\cdot  s^{-1} \]
 have the same eigenvalues.   	
\end{lemma}

\begin{proof}
 The algebra isomorphism $\phi : A \to B$ induces a linear isomorphism $\widetilde{\phi}: I/I^2 \to J/J^2$, 
 $\widetilde{a}  \mapsto \widetilde{\phi(a)}$, with for $a_1,a_2\in A$ and  $x \in I_A$
 \[
  \widetilde{\phi}(\overline{a_1}\cdot \tilde{x}\cdot \overline{a_2}) = \widetilde{\phi}(\widetilde{a_1xa_2})= 
  \widetilde{\phi(a_1xa_2)}= \widetilde{\phi(a_1)\phi(x)\phi(a_2)}= \overline{\phi(a_1)}\cdot \widetilde{\phi(x)}\cdot \overline{\phi(a_2)}
 \]
 Hence for the invertible elements $t = \overline{a}$ and $s=\overline{\phi(a)}$, we have 
 $\widetilde{\phi} \circ {\rm ad}(t)= {\rm ad}(s) \circ \widetilde{\phi}$, and our assertion follows. 
\end{proof}

Lemma \ref{lem:conormal} will be applied with $I$ the commutator ideal of $A$, and the perfect situation will be when the abelianization of $A$ is a one variable Laurent polynomials algebra, or a product of such algebras, because of the following easy well-known lemma.

\bl\label{lem:autolaurent}
\begin{enumerate}
 \item For any algebra automorphism $\varphi : k[t^{\pm 1}] \to k[t^{\pm 1}]$, there exists $\lambda \in k^*$ such that $\varphi(t)= \lambda t^{\pm 1}$.
 \item For any algebra automorphism $\varphi : k[t^{\pm 1}] \times k[t^{\pm 1}] \to k[t^{\pm 1}]\times k[t^{\pm 1}]$, we have $\varphi(t,1) \in \{(\lambda t^{\pm 1}, 1), (1, \lambda t^{\pm 1}), \ \lambda \in k^*\}$.
\end{enumerate}
\el

\bpf
The group of invertible elements in  $k[t^{\pm 1}]$  is $\{\lambda t^n, \lambda \in k^*, n \in \mathbb Z\}$, hence  a generator argument for this  group yields the first claim.  To prove the second one,  first notice that, since an automorphism of an algebra must permute its idempotents, an automorphism $R \times R \to R \times R$ for some ring $R$ with non trivial idempotents is either of the form $(x,y) \mapsto (f_1(x), f_2(y))$ or of the form $(x,y) \mapsto (f_1(y), f_2(x))$ for some algebra automorphisms $f_1, f_2$ of $R$. In particular, by the first claim, an automorphism $\varphi$ of $k[t^{\pm 1}]\times k[t^{\pm 1}]$ satisfies $\varphi(t,1) \in \{(\lambda t^{\pm 1}, 1), (1, \lambda t^{\pm 1}), \ \lambda \in k^*\}$. 
\epf

A typical and simple example of an application of Lemma \ref{lem:conormal} is the following certainly well-known result.

\bex \label{ex:Aq}
For $q \in k^*$, let $A_q=k\langle x^{\pm 1}, \ y  \ | \ yx=qxy\rangle$. For $p,q \in k^*$, we have $A_q\simeq A_p \iff p=q^{\pm 1}$.
\eex

\bpf
For $p=q^{-1}$, we have $A_q\simeq A_{q^{-1}}$, with an isomorphism obtained by sending $x$ to $x^{-1}$ and $y$ to $y$. 

Let $I$ be the commutator ideal of $A_q$. If $q\not = 1$, it is easily seen that $I=Ay$ and that $A/I \simeq k[t^{\pm 1}]$ with $u=\overline{x}$ corresponding to $t$,  while $I/I^2 = (A/I)\widetilde{y}$. We then have, since $A/I$ is commutative and $\widetilde{y}\not = 0$ (as one sees using the usual basis of $A_q$), $\textrm{Sp}({\rm ad}(u))=\{q^{-1}\}$.

Assume now given an algebra isomorphism $\phi :  A_q=A \to B=A_p$. If $q=1$ then $p=1$ because these are the only cases where the algebras are commutative. So assume that $q\not =1$ and $p\not = 1$. By Lemma \ref{lem:autolaurent} the induced isomorphism $\overline{\phi} : A_{\rm ab} \to B_{\rm ab}$ satisfies $\overline{\phi}(u) =\lambda u^{\pm 1}$, for $\lambda \in k^*$, hence ${\rm Sp}({\rm ad}(\overline{\phi}(u)))=\{p\}$ or  ${\rm Sp}({\rm ad}(\overline{\phi}(u)))=\{p^{-1}\}$. We conclude, by Lemma \ref{lem:conormal}, that $q=p$ or $q=p^{-1}$.
\epf

As another application of Lemma \ref{lem:conormal}, we mention the following result, proved in \cite{gad14} (and generalized in \cite{SAQ15}).

\bex \label{ex:qWeyl}
For $q \in k^*$, consider the quantum Weyl algebra \[A_1(q)=k\langle x^{\pm 1}, \ y  \ | \ xy-qyx=1\rangle\] For $p,q \in k^*$, we have $A_1(q)\simeq A_1(p) \iff p=q^{\pm 1}$.
\eex

\bpf
For $p=q^{-1}$, we have $A_1(q)\simeq A_1(q^{-1})$, by \cite[Proposition 6.1]{gad14} (the isomorphism is given by $\phi(x)=q^{-1}y$ and $\phi(y)=-x$).

Let $z= xy-yx \in A= A_1(q)$. We have $xz=qzx$ and $yz=q^{-1}zy$, hence $I=Az=zA$ is the commutator ideal in $A$. Let $u=\overline{x} \in A_{\rm ab}=A/I$. Assume that  $q \not=1$. Then $u$ is invertible, with inverse $(1-q)\overline{y}$. We have $I/I^2 = (A/I) \widetilde{z}$, and  
\[{\rm ad}(u)(\widetilde{z})= (1-q)\overline{x} \cdot \widetilde{z}\cdot \overline{y} =
(1-q) \widetilde{xzy} = (1-q)q\overline{xy}\cdot \widetilde{z}= q\widetilde{z} \]
Since $A/I$ is commutative and $\widetilde{z} \not=0$ (we see that $z\not \in I^2$ using  that $\{z^ix^jy^k, \ i,j,k \in \mathbb N\}$ is a basis of $A_1(q)$, which follows by using the diamond lemma) we get that ${\rm Sp}({\rm ad}(u))=\{q\}$. Since $A/I \simeq k[t^{\pm 1}]$, we now conclude exactly as in Example \ref{ex:Aq}.
%
\epf

\subsection{Proof of Theorem \ref{thm:iso}}
Let $p,q,\alpha, \beta \in k^*$, and assume given an isomorphism 
\[\phi : \OO=A \longrightarrow B = \mathscr{O}_{\alpha, \beta}(\normalfont\text{SL}_{2}(k))\]
Let $I$ be the commutator ideal of $A$, so that $J =\phi(I)$ is the commutator ideal of $B$, and that $\phi$ induces an isomorphism 
\[\overline{\phi} : A_{\rm ab} =A/I \longrightarrow B/J = B_{\rm ab}\] 
We now want to use Lemma \ref{lem:conormal}, and we
discuss according to the dichotomy in Proposition \ref{prop:abelianization}.

\subsubsection{The case $p\not=1$, $q\not=1$ and $pq\not=1$}
It follows from Proposition \ref{prop:abelianization} that $\alpha$, $\beta$ satisfy the same conditions, and  
$A_{\rm ab} \simeq k[t^{\pm 1}]\simeq B_{\rm ab}$. 
By Lemma \ref{lem:autolaurent} we have
$\overline{\phi}(\overline{a}) = \overline{\phi(a)}=\lambda \overline{a}^{\pm 1}$, for some $\lambda \in k^*$.
We easily see, using the basis in Proposition  \ref{prop:diamond2}, that  we have  $I/I^2= A/I\cdot \widetilde{b} + A/I \cdot \widetilde{c}$ as an $A/I$-module, with $\widetilde{b}$ and $\widetilde{c}$ non-zero elements in $I/I^2$. We have
\begin{align*} {\rm ad}(\overline{a})(\widetilde{b}) &= \overline{a} \cdot \widetilde{b} \cdot \overline{a}^{-1}= \overline{a} \cdot \widetilde{b} \cdot \overline{d}= \widetilde{abd} = q^{-1} \widetilde{adb} = q^{-1}\overline{ad}\cdot\widetilde{b}=q^{-1}\widetilde{b} \\
\text{and} \quad   
{\rm ad}(\overline{a})(\widetilde{c}) &= \overline{a} \cdot \widetilde{c} \cdot \overline{a}^{-1}= \overline{a} \cdot \widetilde{c} \cdot \overline{d}= \widetilde{acd} = p^{-1} \widetilde{adc} = p^{-1}\overline{ad}\cdot\widetilde{c}=p^{-1}\widetilde{c}
\end{align*}
Since $A/I$ is commutative, we deduce that ${\rm Sp}({\rm ad}(\overline{a}))=\{p^{-1}, q^{-1}\}$. Similarly, in $B$, the same computation gives ${\rm Sp}({\rm ad}(\overline{\phi(a)}))=\{\alpha^{-1}, \beta^{-1}\}$ or 
${\rm Sp}({\rm ad}(\overline{\phi(a)}))=\{\alpha, \beta\}$. Hence by Lemma \ref{lem:conormal}, we have 
\[\{p^{-1},q^{-1}\} = \{\alpha, \beta\} \ \text{or} \ \{p,q\} = \{\alpha, \beta\},\]
as expected.

\subsubsection{The case $p=1$, $q\not=1$} By Proposition \ref{prop:abelianization} and Proposition \ref{prop:isobasic}
we can assume that $\alpha=1$ and $\beta\not=1 $, and 
$A_{\rm ab} \simeq k[t^{\pm 1}, y]\simeq B_{\rm ab}$.
The group of invertible elements in  $k[t^{\pm 1},y]$   is $\{\lambda t^n, \lambda \in k^*, n \in \mathbb Z\}$, hence   as in Lemma \ref{lem:autolaurent} a generator argument for this group yields that
$\overline{\phi}(\overline{a}) = \overline{\phi(a)}=\lambda \overline{a}^{\pm 1}$, for some  $\lambda \in k^*$. 
We have $I/I^2 = A/I\cdot \widetilde{b}$, and  exactly as in the previous case, we see that ${\rm Sp}({\rm ad}(\overline{a}))=\{q^{-1}\}$, and conclude that $\beta =q^{\pm 1}$.

\subsubsection{The case $q=p^{-1}\not=1$} Again, by Proposition \ref{prop:abelianization} we have $\beta=\alpha^{-1}$, with also
$A_{\rm ab} \simeq k[t^{\pm 1}]\times  k[t^{\pm 1}]\simeq B_{\rm ab}$, and 
$I =Aab +Aac+Abd+Acd$. In $A_{\rm ab}=A/I$, put \[u = \overline{a-p^{-1}bc}\]
It is straightforward that $u$ is invertible with $u^{-1}= \overline{d-p^{-1}bc}$. In order to compute ${\rm Sp}({\rm ad}(u))$ we consider $e=\overline{ad}$ and $f=-p^{-1}\overline{bc}$. These are orthogonal idempotents in $A/I$, with $\overline{1}=e+ f$. We thus have 
\begin{align*}I/I^2 &= (A/I)e\cdot\widetilde{ab} + (A/I) f\cdot \widetilde{ab} + (A/I) e\cdot\widetilde{ac}+ (A/I) f\cdot \widetilde{ac}
\\
&+(A/I)e\cdot \widetilde{bd} + (A/I)f\cdot \widetilde{bd}+(A/I)e\cdot \widetilde{cd} +(A/I) f\cdot \widetilde{cd}\end{align*}
with the $8$ elements $e\cdot \widetilde{ab}, \ldots , f \cdot \widetilde{cd}$ all non zero (by Proposition \ref{prop:diamond2}).
We have
\[ue= \overline{a} =eu \quad  \textrm{and} \quad uf = f = fu
\]
and hence 
\[  {\rm ad}(u)(e\cdot \widetilde{ab}) = ue\cdot \widetilde{ab} \cdot u^{-1} = \widetilde{a^2b}\cdot u^{-1} = p \widetilde{ab} \cdot \overline{a}u^{-1} =
p \widetilde{ab} \cdot e = p e \cdot \widetilde{ab}
\] 
Similar computations give
\[{\rm ad}(u)(f\cdot \widetilde{ab})= f\cdot \widetilde{ab}, \
{\rm ad}(u)(e\cdot \widetilde{ac})= p^{-1}e\cdot \widetilde{ac}, \ 
{\rm ad}(u)(f\cdot \widetilde{ac})= f\cdot \widetilde{ac},
\]
\[{\rm ad}(u)(e\cdot \widetilde{bd})= pe\cdot \widetilde{bd}, \
{\rm ad}(u)(f\cdot \widetilde{bd})= f\cdot \widetilde{bd}, \
{\rm ad}(u)(e\cdot \widetilde{cd})= p^{-1}e\cdot \widetilde{cd},
{\rm ad}(u)(f\cdot \widetilde{cd})= f\cdot \widetilde{cd}\]
and hence (again using that $A/I$ is commutative) we obtain ${\rm Sp}({\rm ad}(u))=\{1, p, p^{-1}\}$. Define now 
\[v = \overline{ad+b} \in A/I\]
This is an invertible element, with $v^{-1}=\overline{ad-p^{-1}c}$.
We have $ve=e=ev$ and $vf=\overline{b}=fv$, and similarly to what we did for $u$, we check that ${\rm Sp}({\rm ad}(v))=\{1, p, p^{-1}\}$. 

Since under the natural isomorphism  $A_{\rm ab}\simeq k[t^{\pm 1}]\times k[t^{\pm 1}]$ in Proposition \ref{prop:abelianization}, the elements $u$ and $v$ correspond respectively to $(t,1)$ and $(1,t)$, by Lemma \ref{lem:autolaurent}, the isomorphism $\overline{\phi} : A_{\rm ab} \simeq B_{\rm ab}$ satisfies $\overline{\phi}(u)\in \{\lambda u^{\pm 1}, \lambda v^{\pm 1}, \lambda \in k^*\}$. Thanks to Lemma \ref{lem:conormal}, this yields $\{1, p, p^{-1}\} = \{1, \beta, \beta^{-1}\}$, as we expect, and this concludes the proof of Theorem \ref{thm:iso}.

\brq
Lemma \ref{lem:conormal} can be used as well to recover D'Andrea's isomorphism   result  \cite{Dan24} for quantum spheres, since in this case again the abelianization is an algebra of Laurent polynomials in one variable, and the deformation parameter is the eigenvalue of the corresponding generator.
\erq

\brq
In lemma \ref{lem:conormal} we considered an invertible element $t=\overline{a}$ in $A/I$  and the automorphism $\mathrm{ad}(t)$ of $I/I^2$. It is of course possible to write a version of Lemma \ref{lem:conormal} without assuming this,  which will be more general but less pleasant, we feel. This can be formulated as follows (we retain all the notation in Lemma \ref{lem:conormal}). For $a \in A$, denote by $R_{\overline{a}}$ and $L_{\overline{a}}$ the respective endomorphisms of $I/I^2$ defined by $L_{\overline{a}}(\widetilde{x})=\widetilde{ax}$ and $R_{\overline{a}}(\widetilde{x}) = \widetilde{xa}$. The statement in Lemma \ref{lem:conormal} becomes
\[\left\{\lambda \in k^*, \ \mathrm{Ker}(L_{\overline{a}} - \lambda R_{\overline{a}}) \not=\{0\}\right\} 
= \left\{\lambda \in k^*, \ \mathrm{Ker}(L_{\overline{\phi(a)}} - \lambda R_{\overline{\phi(a)}}) \not=\{0\}\right\} \]
\erq

\bibliographystyle{plain}
\bibliography{biblio_iso.bib}

\end{document}